\documentclass[a4paper,11pt]{article}

\usepackage{graphicx} 
\usepackage{amsmath}
\usepackage{amssymb}
\usepackage{amsfonts}
\usepackage{hyperref} 
\usepackage{amsthm}
\usepackage{mathtools}

\usepackage{tikz}
\usetikzlibrary{graphs,graphs.standard}

\usepackage[utf8]{inputenc}
\usepackage{graphicx}

\usepackage{float}
\usepackage{caption}
\usepackage{diagbox} 

\usepackage[verbose]{placeins}

\usepackage[numbers]{natbib}
\usepackage{cite}

\usepackage{setspace}

\usepackage{multirow}

\theoremstyle{plain} 

\theoremstyle{definition} 
\newtheorem{defn}{Definition}[section]
\newtheorem{eg}{Example}[section]

\theoremstyle{remark} 
\newtheorem{rem}{Remark}[section]

\title{Maximal cross-ratio degree for 8 points in $\mathbb{P}^1$}
\author{Arjun Maniyar}
\date{}

\begin{document}

\maketitle

\begin{abstract}
     The cross-ratio degree problem asks for the number of configurations of $n$ points in $\mathbb{P}^1$ that satisfy $n-3$ specified cross-ratio conditions. Silversmith points out several connections and applications of cross-ratio degree in \citep{silversmith2022perfectmathing,silversmith2024triangulation}. It is also known that the maximal cross-ratio degree for 8 points is at least 4, see \citep{silversmith2024triangulation}. In this paper, we will see that the maximal cross-ratio degree for 8 points in $\mathbb{P}^1$ is equal to 4. 
\end{abstract}

\section{Introduction}\label{section1}

We will use the notation and results from \citep{silversmith2022perfectmathing}  throughout this paper. 
\begin{defn}
    Let $z_1,z_2,z_3,z_4$ be distinct points in the projective space $\mathbb{P}^1$, then their \textit{cross-ratio} is defined by 
\[
\mathrm{CR}(z_1,z_2,z_3,z_4)
= \frac{(z_3 - z_1)(z_4 - z_2)}{(z_3 - z_2)(z_4 - z_1)} \in \mathbb{C}\setminus\{0,1\}.
\]
\end{defn}

Informally, cross-ratio degree is the number of configurations of points in $\mathbb{P}^1$ given some constraints on their cross-ratios. This can be formalised as follows \citep[p.~1]{silversmith2022perfectmathing}. 
Let $n \ge 4$. Fix a collection
\[
\mathcal{T} = (e_1,\dots,e_{n-3}),
\]
where each $e_i = (v_{i,1},v_{i,2},v_{i,3},v_{i,4})$ is a 4-tuple of distinct elements of $\{1,\dots,n\}$.
For any configuration of distinct points $(z_1,\dots,z_n) \in (\mathbb P^1)^n$, define
\[
\begin{aligned}
\mathrm{CR}_{\mathcal{T}}(z_1,\ldots,z_n)
=
\bigl(
&\mathrm{CR}(z_{v_{1,1}},z_{v_{1,2}},z_{v_{1,3}},z_{v_{1,4}}),
\\
&\mathrm{CR}(z_{v_{2,1}},z_{v_{2,2}},z_{v_{2,3}},z_{v_{2,4}}),
\\
&\qquad\vdots
\\
&\mathrm{CR}(z_{v_{n-3,1}},z_{v_{n-3,2}},z_{v_{n-3,3}},z_{v_{n-3,4}})
\bigr)
\in (\mathbb{C}\setminus\{0,1\})^{\,n-3}.
\end{aligned}
\]
We define the the space $M_{0,n}$ as in \citep{silversmith2022perfectmathing,silversmith2024triangulation}. Define $M_{0,n}$ to be the the moduli space of $n$ ordered, distinct points
$(z_1,\dots,z_n)\in (\mathbb P^1)^n$, up to M\"obius transformations. Note that $M_{0,n}$ is a $(n-3)$-dimensional variety. Since cross-ratios are invariant under M\"obius transformations,
$\mathrm{CR}_{\mathcal{T}}$ defined above gives us a map from $M_{0,n} \to (\mathbb{C}\setminus\{0,1\})^{n-3}$. 
\begin{defn}\citep[p.~1]{silversmith2022perfectmathing}
    The \emph{cross-ratio degree} of $\mathcal{T}$, denoted $d_\mathcal{T}$, is defined to be the degree of  $\mathrm{CR}_{\mathcal{T}} : M_{0,n} \to (\mathbb{C}\setminus\{0,1\})^{n-3}$. This is well-defined since $\dim(M_{0,n}) = \dim((\mathbb{C}\setminus\{0,1\})^{n-3}) = n-3$.
\end{defn}

Observe that the following operations does not affect $d_{\mathcal{T}}$ \citep[p.~1]{silversmith2022perfectmathing}:
\begin{itemize}
    \item permuting the 4 elements in $e_i$,
    \item permuting the sets $e_i$  in $\mathcal{T} = (e_1, \dots, e_{n-3})$,
    \item renaming the elements of the set $\{1, \dots , n\}$.
\end{itemize}
Now think of the points in $\mathbb{P}^1$ as the vertices of a hypergraph and each cross-ratio relation as a hyperedge, which gives us a hypergraph $\mathcal{T} = (V,E)$ where $|V| = n$, $|E| = n-3$, and $|e| =4$ for every $e \in E$. From the observations we can deduce that the cross-ratio degree $d_\mathcal{T}$ depends only on the isomorphism class of the hypergraph $\mathcal{T} = (V,E)$. 

The problem at hand is to calculate the maximal cross-ratio degree for 8 points in $\mathbb{P}^1$. One way is to calculate cross-ratio degree for all non isomorphic hypergraphs and find the maximum. Let $s_{p,k}^n$ be the number of (non isomorphic) $(n+1)$-uniform hypergraphs with $p$ vertices and $k$ hyperedges. Note that in our case every vertex of the hypergraph must be part of at least one edge, but $s_{p,k}^n$ also counts graphs with isolated vertices. Thus, the number of hypergraphs with our requirements on 8 vertices is given by $s_{8,5}^3 - s_{7,5}^3$. See Section \ref{number of hypergraphs}(Appendix) on how these numbers can be calculated using results from \citep{PALMER}. These numbers are 621 and 137 respectively, so, the number of hypergraphs is $621-137 = 484$. 

In \citep[p.~2]{silversmith2024triangulation}, Silversmith provides lower bounds on the maximal cross-ratio degrees for up to 14 points in $\mathbb{P}^1$. For the case of 8 points, the bound provided is $\geq 4$. We will see that this bound is exact, that is, the maximal cross-ratio degree is 4. 
\begin{rem} \label{remark:7 points cross-ratio degree}
    For 7 points, the lower bound provided in \citep[p.~2]{silversmith2024triangulation} is $\geq 2$ which is also exact, see \citep{Powell2024}. And, the number of valid hypergraphs on 7 vertices is $s_{7,4}^3 - s_{6,4}^3 = 29$. 
\end{rem}

\section{Calculating cross-ratio degree}

There are several ways of calculating cross-ratio degree, Silversmith outlines six different ways in \citep{silversmith2022perfectmathing}. One of these methods is the naive approach of solving the cross-ratio equations directly and counting the number of solutions. Since we are only concerned with small cases, that is 8 points with 5 cross-ratio equations, we shall only use this simple and naive approach. This method can be best explained through an example. 

\begin{eg}
    Let $\mathcal{T} = (V,E)$ with $V = \{1,\ldots,8\}$ and 
$E = \{\{1,2,3,4\},\\ \{1,2,6,7\}, \{1,3,7,8\}\, \{1,2,5,8\}\, \{3,4,5,6\}\}$. Let  $G_{\mathcal{T}} = (V \sqcup E, I)$ be the bipartite graph with parts $V$ and $E$, and whose edges are given by the set $I = \{(v,e) \in V \times E : v \in e\}$. Then the biadjacency matrix of $G_{\mathcal{T}}$ is given by (where the top row is the column sum)
\[ 
\begin{pmatrix}
4 & 3 & 3 & 2 & 2 & 2 & 2 & 2 \\
1 & 1 & 1 & 1 & 0 & 0 & 0 & 0 \\
1 & 1 & 0 & 0 & 0 & 1 & 1 & 0 \\
1 & 0 & 1 & 0 & 0 & 0 & 1 & 1 \\
1 & 1 & 0 & 0 & 1 & 0 & 0 & 1 \\
0 & 0 & 1 & 1 & 1 & 1 & 0 & 0
\end{pmatrix}.
\]
We will calculate $d_{\mathcal{T}}$ by counting the number of tuples $(p_1,\ldots,p_8) $ satisfying
\[
\begin{gathered}
\mathrm{CR}(p_1,p_2,p_3,p_4) = a_1 \qquad
\mathrm{CR}(p_1,p_2,p_6,p_7) = a_2 \qquad
\mathrm{CR}(p_1,p_3,p_7,p_8) = a_3
\\
\mathrm{CR}(p_1,p_2,p_5,p_8) = a_4 \qquad
\mathrm{CR}(p_3,p_4,p_5,p_6) = a_5,
\end{gathered}
\]
where $a_1,a_2,a_3,a_4,a_5$ are generically chosen scalars. Since cross ratio is invariant under projective transformations, we can change coordinates so that $p_1 = \infty$, $p_2 = 0$, $p_3 = 1$.
Then the equations above simplify to
\[
p_4 = a_1
\qquad
\frac{p_7}{p_6} = a_2
\qquad
\frac{p_8-1}{p_7-1} = a_3
\qquad
\frac{p_8}{p_5} = a_4
\qquad
\frac{(p_5 - 1)(p_6 - p_4)}{(p_6 - 1)(p_5 - p_4)} = a_5.
\]
Note that the first four equations are linear, so, we can write $p_5$ linearly in terms of $p_6$: $$p_5=\frac{(a_2p_6-1)a_3+1}{a_4}.$$
Then substituting this equation and $p_4=a_1$ in the last equation yields a quadratic equation in terms of $p_6$. This gives two distinct tuples which satisfy the equations above. So, the total number of solutions for this system is 2.
Therefore, $d_{\mathcal{T}} = 2$.
\end{eg}

\begin{rem}
    
    As in the example above, every hypergraph can be represented in form of biadjacency matrix. So, from this point on, we will use matrix representations for hypergraphs, without stating it explicitly each time.
\end{rem}

We now make some observations.\\
\textbf{Observation 1:} Swapping columns in the biadjacency
matrix is equivalent to renaming the elements of $V$, which does not affect the cross-ratio degree. Thus, we may reorder the columns so that their sums are nonincreasing
from left to right (we already had the matrix in this form in the example above). In particular, the column with the largest sum can be
placed first, and we can set $p_1 = \infty$. Then any cross-ratio relation involving
$p_1$ then becomes linear, since the $\infty$ terms cancel in numerator and
denominator. This means that if any matrix  has the column
\[ 
\begin{pmatrix}
1 \\
1  \\
1  \\
1  \\
1  

\end{pmatrix},
\]
then that system of equations has at most one solution since all the equations are linear. So, any matrix which has a column of sum 5 has $d_{\mathcal{T}} \leq 1$. This idea is originally from \citep[p.~15]{Powell2024}. Similarly, if the largest column sum is $4$ (i.e.\ a column with four $1$'s and one $0$), then we obtain a system of four linear equations and at most one quadratic equation. In this case, we therefore have $d_{\mathcal{T}} \le 2$ (as in the previous example). However, we cannot take this idea any further: with largest column sum $3$ one can obtain cubic equations or two quadratics; see Section~\ref{appendix:data} for such matrices.

All possible column sums for the biadjacency matrices are listed in
Table~\ref{table:number of matrices for col sum} in
Section~\ref{appendix:data}. Since we have already dealt with column sums
containing $5$ and $4$, the only remaining cases of interest are those with
the following column sums:
\[
(3,3,3,3,2,2,2,2),\quad
(3,3,3,3,3,2,2,1),\quad
(3,3,3,3,3,3,1,1).
\]

\textbf{Observation 2:}
Recall that \textit{degree} of a vertex in a hypergraph is defined to be the number of edges that contain that vertex. Now consider a hypergraph that has at least one vertex of degree 1, but does not contain an edge that contains more than one vertex of degree 1 (see Remark \ref{rem: edge with multiple degree 1} for why this condition is needed).
So, in the matrix we have a column with one 1 and four 0's (column sum is 1). Note that each row in the matrix corresponds to a hyperedge in the hypergraph, so, swapping rows of a matrix does not affect the cross-ratio degree since it does not affect the cross-ratio constraints. Thus, we may swap rows to get the matrix in the following form:  
\[ M=
\left(
\begin{array}{c|c}
A & \begin{array}{c} 0 \\ 0 \\ 0 \\ 0 \end{array}
\\ \hline
B & 1
\end{array}
\right),
\]
where $A$ is a $4 \times 7 $ matrix and $B$ is a $1 \times 7$ matrix. Now delete the last column and the last row from the matrix $M$, which leaves us with the matrix $A$. Observe that the matrix $A$ represents a 4-uniform hypergraph on 7 vertices with no isolated vertices. So, the matrix $A$ leaves us with a ``valid" cross-ratio degree problem for 7 points. And, the largest cross-ratio degree for 7 points is 2 (see Remark \ref{remark:7 points cross-ratio degree}). This means that there are at most 2 tuples $(p_1,\dots ,p_7)$ satisfying the cross-ratio constraints. Now we add back the last column and row. From \citep[Lemma~3.2]{silversmith2022perfectmathing}, it follows that $d_{\mathcal{T}} = d_{\mathcal{T}'}$ where $\mathcal{T}$ and $\mathcal{T}'$ are hypergraphs corresponding to matrices $M$ and $A$ respectively. Informally put, we can uniquely determine $p_8$ since we can write $p_8$ in terms of $p_1,\dots ,p_7$ using the cross-ratio equation from the last row. Thus, there are at most 2 solutions to the system given by the matrix $M$. 

\begin{rem}\label{rem: edge with multiple degree 1}
     Note that the argument above does not work if we have an edge that contains more than one vertex of degree 1. This is because we get the matrix in the following form:
\[ M=
\left(
\begin{array}{c|c c}
A & 
\begin{array}{c} 0 \\ 0 \\ 0 \\ 0 \end{array} &
\begin{array}{c} 0 \\ 0 \\ 0 \\ 0 \end{array}
\\ \hline
B & 1 & 1
\end{array}
\right).
\]
     Deleting the last row and column in the above matrix $M$ does not reduce the problem to a ``valid" cross-ratio degree matrix for 7 points since the last column is 0. In fact, in these cases we get $d_{\mathcal{T}} = 0$ by \citep[Theorem~1.1]{silversmith2022perfectmathing}, because after deleting $p_1,p_2,p_3$ there are no perfect matchings. Another way to see that $d_{\mathcal{T}} = 0$ is to note that the first four rows of $M$ yield a relation among the $a_i$'s. But since $a_1,\dots, a_5$ are chosen generically, we get $d_{\mathcal{T}} = 0$. 
\end{rem}

We know that the matrices classified in Remark \ref{rem: edge with multiple degree 1} have degree 0. But these are not the only matrices with degree 0, for example consider the following
\[
\begin{pmatrix}
3 & 3 & 3 & 3 & 2 & 2 & 2 & 2 \\
1 & 1 & 1 & 1 & 0 & 0 & 0 & 0 \\
1 & 1 & 1 & 0 & 1 & 0 & 0 & 0 \\
1 & 1 & 0 & 1 & 1 & 0 & 0 & 0 \\
0 & 0 & 1 & 0 & 0 & 1 & 1 & 1 \\
0 & 0 & 0 & 1 & 0 & 1 & 1 & 1
\end{pmatrix}.
\]
As usual we set $p_1 = \infty$, $p_2 = 0$, $p_3 = 1$. Then from the first three rows we get the following constraints:  
\[
p_4 = a_1
\qquad
p_5 = a_2
\qquad
\frac{p_5}{p_4} = a_3.
\]
Thus, we must have $\frac{a_2}{a_1} = a_3$. But since the values $a_1,\dots ,a_5$ are chosen generically, we must have that $d_{\mathcal{T}} = 0$. This can also be shown using \citep[Theorem~1.1]{silversmith2022perfectmathing} as before to get $d_{\mathcal{T}} \leq 0$. Note that all these cases of $d_{\mathcal{T}} = 0$ fall under \citep[Proposition~4.1]{silversmith2022perfectmathing}. 

The observations mentioned above are not restricted to 8 points but are also true for arbitrary number of points in $\mathbb{P}^1$. Finally, from these observations we have that all matrices, except those with column sum $(3,3,3,3,2,2,2,2)$, have cross-ratio degree at most 2. Computing the cross-ratio degree of all matrices with column sum $(3,3,3,3,2,2,2,2)$ yields that the maximal cross-ratio degree is 4. It is also worth noting that there are matrices with column sum $(3,3,3,3,2,2,2,2)$ with $d_{\mathcal{T}} = 0,1,2$  and 3. Matrices with degree 3 and 4 are listed in Section \ref{appendix:data}.

\subsection{Computing all biadjacency matrices}

To find all 484 non-isomorphic hypergraphs, we can use the notion of a graph canonical form (any two isomorphic hypergraphs have the same canonical form). A widely known canonical form is defined as the lexicographically smallest representative of an isomorphism class, namely the graph whose adjacency matrix, viewed as a linear string, is lexicographically smallest among all graphs in the class. It is straightforward to define and compute, however, it is not computationally efficient. By generating hypergraphs via backtracking and pruning methods, and applying an isomorphism check based on canonical forms, we can compute all required hypergraphs. While this approach is not optimal in general, it is sufficient for the small cases considered in this paper, which is hypergraphs on 8 vertices.

\section{Acknowledgements} 
This work was undertaken as part of my URSS (Undergraduate Research Support Scheme) funded summer research project at the Warwick Mathematics Institute, carried out under the supervision of Helena Verrill. I sincerely thank Helena Verrill for their guidance and feedback throughout the project.
\section{Appendix}

\subsection{Counting number of $n$-uniform hypergraphs}\label{number of hypergraphs}

We first define an $n$-plex as follows \citep[p.~1]{PALMER}.

\begin{defn}
    An \textit{$n$-plex} of order $p$ is an $n$-dimensional simplicial complex on $p$ vertices such that each maximal simplex has dimension either $n$ or $0$. (A maximal simplex is a simplex that is not a face of any larger simplex).
\end{defn}

For instance, a graph is a 1-plex. Palmer provides a method to enumerate these objects in \citep{PALMER}. Let $s_p^n(x)$ be the counting polynomial for $n$-plexes of order $p$, so we have \[s_p^n(x)=\sum s_{p,k}^n x^k,\] where $s_{p,k}^n$ is the number of $n$-plexes with $k$ simplexes of dimension $n$. Now note that $s_{p,k}^n=$ the number of (non isomorphic) $(n+1)$-uniform hypergraphs with $p$ vertices and $k$ hyperedges. It is also known that (see \citep[pp.~25, 32]{harary1967} for a proof) 
$$s_p^n(x) = Z(S_p^{(n+1)}, 1+x),$$ 
where $Z$ is the cycle index and $S_p^{(n+1)}$ is an induced permutation group (see \citep{PALMER} for more details).
Using the method provided to calculate the cycle index of an induced permutation group in \citep{PALMER}, we can calculate $s_p^n(x)$ and hence $s_{p,k}^n$.
\subsection{Matrices and other data} \label{appendix:data}
This section contains all the computational data that might be worth sharing.
Recall that the total number of valid biadjacency matrices is 484. As before, the top row in the following matrices is dedicated to the column sum.\\

Matrices with $d_{\mathcal{T}} = 4$: 
\[
\begin{pmatrix}
3 & 3 & 3 & 3 & 2 & 2 & 2 & 2 \\
1 & 1 & 1 & 1 & 0 & 0 & 0 & 0 \\
1 & 1 & 0 & 0 & 1 & 1 & 0 & 0 \\
1 & 0 & 1 & 0 & 0 & 0 & 1 & 1 \\
0 & 1 & 0 & 1 & 0 & 0 & 1 & 1 \\
0 & 0 & 1 & 1 & 1 & 1 & 0 & 0
\end{pmatrix}
\quad
\begin{pmatrix}
3 & 3 & 3 & 3 & 2 & 2 & 2 & 2 \\
1 & 1 & 1 & 0 & 1 & 0 & 0 & 0 \\
1 & 1 & 0 & 1 & 0 & 1 & 0 & 0 \\
1 & 0 & 1 & 0 & 0 & 0 & 1 & 1 \\
0 & 1 & 0 & 1 & 0 & 0 & 1 & 1 \\
0 & 0 & 1 & 1 & 1 & 1 & 0 & 0
\end{pmatrix}
\]
\[
\begin{pmatrix}
3 & 3 & 3 & 3 & 2 & 2 & 2 & 2 \\
1 & 1 & 1 & 1 & 0 & 0 & 0 & 0 \\
1 & 1 & 0 & 0 & 1 & 1 & 0 & 0 \\
1 & 1 & 0 & 0 & 0 & 0 & 1 & 1 \\
0 & 0 & 1 & 1 & 1 & 1 & 0 & 0 \\
0 & 0 & 1 & 1 & 0 & 0 & 1 & 1
\end{pmatrix}
\quad
\begin{pmatrix}
3 & 3 & 3 & 3 & 2 & 2 & 2 & 2 \\
1 & 1 & 1 & 0 & 1 & 0 & 0 & 0 \\
1 & 1 & 0 & 1 & 0 & 1 & 0 & 0 \\
1 & 1 & 0 & 0 & 0 & 0 & 1 & 1 \\
0 & 0 & 1 & 1 & 1 & 1 & 0 & 0 \\
0 & 0 & 1 & 1 & 0 & 0 & 1 & 1
\end{pmatrix}
\]
\\

Matrices with $d_{\mathcal{T}} = 3$:

\[
\begin{pmatrix}
3 & 3 & 3 & 3 & 2 & 2 & 2 & 2 \\
1 & 1 & 1 & 0 & 1 & 0 & 0 & 0 \\
1 & 1 & 0 & 1 & 0 & 1 & 0 & 0 \\
1 & 0 & 1 & 0 & 0 & 1 & 1 & 0 \\
0 & 1 & 0 & 1 & 0 & 0 & 1 & 1 \\
0 & 0 & 1 & 1 & 1 & 0 & 0 & 1
\end{pmatrix}
\quad
\begin{pmatrix}
3 & 3 & 3 & 3 & 2 & 2 & 2 & 2 \\
1 & 1 & 1 & 0 & 1 & 0 & 0 & 0 \\
1 & 1 & 0 & 1 & 0 & 1 & 0 & 0 \\
1 & 0 & 1 & 0 & 0 & 1 & 1 & 0 \\
0 & 1 & 0 & 1 & 1 & 0 & 0 & 1 \\
0 & 0 & 1 & 1 & 0 & 0 & 1 & 1
\end{pmatrix}
\]
\[
\begin{pmatrix}
3 & 3 & 3 & 3 & 2 & 2 & 2 & 2 \\
1 & 1 & 1 & 0 & 1 & 0 & 0 & 0 \\
1 & 1 & 0 & 1 & 0 & 1 & 0 & 0 \\
1 & 1 & 0 & 0 & 0 & 0 & 1 & 1 \\
0 & 0 & 1 & 1 & 1 & 0 & 1 & 0 \\
0 & 0 & 1 & 1 & 0 & 1 & 0 & 1
\end{pmatrix}
\quad
\begin{pmatrix}
3 & 3 & 3 & 3 & 2 & 2 & 2 & 2 \\
1 & 1 & 1 & 0 & 1 & 0 & 0 & 0 \\
1 & 1 & 0 & 1 & 0 & 1 & 0 & 0 \\
1 & 0 & 1 & 1 & 0 & 0 & 1 & 0 \\
0 & 1 & 1 & 1 & 0 & 0 & 0 & 1 \\
0 & 0 & 0 & 0 & 1 & 1 & 1 & 1
\end{pmatrix}
\]
\[
\begin{pmatrix}
3 & 3 & 3 & 3 & 2 & 2 & 2 & 2 \\
1 & 1 & 1 & 0 & 1 & 0 & 0 & 0 \\
1 & 1 & 0 & 1 & 0 & 1 & 0 & 0 \\
1 & 0 & 1 & 1 & 0 & 0 & 1 & 0 \\
0 & 1 & 1 & 0 & 0 & 1 & 0 & 1 \\
0 & 0 & 0 & 1 & 1 & 0 & 1 & 1
\end{pmatrix}
\quad
\begin{pmatrix}
3 & 3 & 3 & 3 & 2 & 2 & 2 & 2 \\
1 & 1 & 1 & 0 & 1 & 0 & 0 & 0 \\
1 & 1 & 1 & 0 & 0 & 1 & 0 & 0 \\
1 & 1 & 0 & 1 & 0 & 0 & 1 & 0 \\
0 & 0 & 1 & 1 & 1 & 0 & 0 & 1 \\
0 & 0 & 0 & 1 & 0 & 1 & 1 & 1
\end{pmatrix}
\]
\[
\begin{pmatrix}
3 & 3 & 3 & 3 & 2 & 2 & 2 & 2 \\
1 & 1 & 1 & 0 & 1 & 0 & 0 & 0 \\
1 & 1 & 1 & 0 & 0 & 1 & 0 & 0 \\
1 & 1 & 0 & 1 & 0 & 0 & 1 & 0 \\
0 & 0 & 1 & 1 & 0 & 0 & 1 & 1 \\
0 & 0 & 0 & 1 & 1 & 1 & 0 & 1
\end{pmatrix}
\quad
\begin{pmatrix}
3 & 3 & 3 & 3 & 2 & 2 & 2 & 2 \\
1 & 1 & 1 & 0 & 1 & 0 & 0 & 0 \\
1 & 1 & 1 & 0 & 0 & 1 & 0 & 0 \\
1 & 0 & 0 & 1 & 1 & 0 & 1 & 0 \\
0 & 1 & 0 & 1 & 0 & 1 & 0 & 1 \\
0 & 0 & 1 & 1 & 0 & 0 & 1 & 1
\end{pmatrix}
\]

\begin{table}[h]
\centering
\caption{Number of matrices for each possible column sum}
\label{table:number of matrices for col sum}

\begin{minipage}{0.48\textwidth}
\centering
\begin{tabular}{|c|c|}
\hline
\multirow{2}{*}{Column sum} & Number \\[-2pt]
& of matrices \\
\hline

$(4,3,3,3,2,2,2,1)$ & 86 \\
$(4,4,3,3,2,2,1,1)$ & 62 \\
$(4,3,3,3,3,2,1,1)$ & 42 \\
$(4,4,3,2,2,2,2,1)$ & 42 \\
$(3,3,3,3,2,2,2,2)$ & 38 \\
$(3,3,3,3,3,2,2,1)$ & 34 \\
$(4,3,3,2,2,2,2,2)$ & 34 \\
$(5,3,3,2,2,2,2,1)$ & 19 \\
$(5,3,3,3,2,2,1,1)$ & 19 \\
$(5,4,3,2,2,2,1,1)$ & 18 \\
$(4,4,3,3,3,1,1,1)$ & 13 \\
$(4,4,4,2,2,2,1,1)$ & 13 \\
$(4,4,4,3,2,1,1,1)$ & 9 \\
$(5,4,3,3,2,1,1,1)$ & 9 \\
\hline
\end{tabular}
\end{minipage}
\hfill
\begin{minipage}{0.48\textwidth}
\centering
\begin{tabular}{|c|c|}
\hline
\multirow{2}{*}{Column sum} & Number \\[-2pt]
& of matrices \\
\hline

$(3,3,3,3,3,3,1,1)$ & 8 \\
$(4,4,2,2,2,2,2,2)$ & 8 \\
$(5,3,2,2,2,2,2,2)$ & 5 \\
$(5,4,2,2,2,2,2,1)$ & 5 \\
$(5,4,4,2,2,1,1,1)$ & 5 \\
$(5,3,3,3,3,1,1,1)$ & 4 \\
$(5,5,2,2,2,2,1,1)$ & 3 \\
$(5,5,3,2,2,1,1,1)$ & 3 \\
$(4,4,4,4,1,1,1,1)$ & 1 \\
$(5,4,4,3,1,1,1,1)$ & 1 \\
$(5,5,3,3,1,1,1,1)$ & 1 \\
$(5,5,4,2,1,1,1,1)$ & 1 \\
$(5,5,5,1,1,1,1,1)$ & 1 \\
\hline
\end{tabular}
\end{minipage}
\end{table}

\begin{table}[h]
\centering
\caption{Distribution of matrices for each $d_{\mathcal{T}} = 1,2,3$ and 4.}
\begin{tabular}{|c | c|}
\hline
$i$ & Number of matrices with $d_{\mathcal{T}} = i$ \\
\hline
$0$ & $79$ \\
$1$ & $279$ \\
$2$ & $114$ \\
$3$ & $8$ \\
$4$ & $4$ \\
\hline
\end{tabular}

\end{table}

\bibliography{references}
\end{document}